\documentclass{elsart}
\usepackage{latexsym,amscd,amsmath}
\usepackage[all]{xy}

\newfont{\pr}{msbm10}
\newfont{\qr}{msbm8}
\newfont{\got}{eufm10}
\renewcommand{\P}{\hbox{\pr P}}
\newcommand{\Z}{\hbox{\pr Z}}
\newcommand{\p}{\hbox{\qr P}}
\newcommand{\I}{\mathcal{I}}

\renewcommand{\O}{\mathcal{ O}}

\newcommand{\ie}{\textrm{i.e.}}

\begin{document}

\begin{frontmatter}

\title{On trisecant lines to White surfaces}
\author{Marie-Am\' elie Bertin}
\address{
12 rue des Bretons,\\
F-94700 Maisons-Alfort (France) \\
\texttt{marie.bertin@parisfree.com}
}

\begin{abstract} Inspired by an argument of Gambier, we show that the only White surface  of ${\P}^{5}$
having a $4$-dimensional trisecant locus is the Segre polygonal surface.
This allows us to deduce that the generic point of the  principal component of
the subvariety $W_{18}[5]$ of  $18$-tuples special in degree $5$ of the  Hilbert scheme of $18$ points
of the plane corresponds  to a smooth $18$-tuple of points in uniform position, not lying on any quartic.
This refines, in this particular case, the general bound due to Coppo.
We also give the number of trisecant lines, counted with multiplicity, which pass through the generic point of a White surface
of non-Segre type.

\end{abstract}

\end{frontmatter}

\section{Introduction}

Let $S$ be a smooth, non-degenerate surface of ${\P}^5$. In order to obtain properties of $S$,
it is classical to study the generic projection of $S$ to ${\P}^3$, 
i.e., the image $\overline{S}$ of $S$ under the projection
from a generic line $L$ of ${\P}^5$. Let $x$ be a point of $\overline{S}$, we denote by $\pi^{-1}(x)$ the scheme theoretical fiber of $\pi$ over $x$. 
Since $\pi$ is a generic projection, $\pi^{-1}$ is $0$-dimensional.
So, the length of the fibers of $\pi$ provides a natural stratification of $\overline{S}$:
\[
\cdots \overline{S}_{k+1}\subseteq \overline{S}_{k}\subseteq \cdots \overline{S}_{1}=\overline{S},
\]
where, for all integers $i\geq 1$,
\[
\overline{S}_{k}:=\{ x\in \overline{S}\, |\, length (\pi^{-1}(x))\geq k\}
\]
is the set of $k$-tuple points of the projection $\pi$.
The following well-known result gives a stratification theorem for surfaces: 
the hypersurface $\overline{S}$ contains a curve of double points, a finite number of triple points
and no quadruple points (\cite{GH} p. 611-618  for this and further properties of the singularities of $\overline{S}$).
Since the projection is generic, over each triple point of $\overline{S}$
lie three distinct points of $S$. Unfortunately, the classical proof of this result, as for instance
in Griffiths and Harris' book \cite{GH} p. 611-618 (or \cite{Be} chapter 9 in arbitrary dimension), is false.
Indeed, as a corollary of the proof those three points
must be in general position (see Dobler \cite{D}). Equivalently, the dimension of the trisecant lines locus of $S$ is at most $3$. 
From the classical proof of the stratification theorem,
one can also deduce the following corollary (\cite{GH} top of page 613):
  `` among $3$-planes of ${\P}^5$ meeting $S$ at points not in general position, the generic one contains $4$ points of $S$ spanning a $2$-plane,
  and not three colinear points''.

Since there exist surfaces of ${\P}^5$ having a $4$-dimensional trisecant line locus, e.g. smooth special Enriques surfaces
(Conte and Verra \cite{CV}, Dolgatchev and Reider \cite{DR}) or Segre polygonal surfaces (Segre \cite{S}, Dobler \cite{D}),
the classical proof of the stratification theorem contains a mistake; we refer to Dobler's thesis \cite{D} for a detailed discussion of this matter.
Following Dobler, we wish to point out, that surfaces of ${\P}^5$ with a $4$-dimensional trisecant lines locus give 
counterexamples to both corollaries of the classical proof of the stratification theorem for surfaces.
The stratification theorem holds nonetheless; indeed, the results of Mather (\cite{Ma1} and \cite{Ma2})
provide a stratification theorem for generic projections to hypersurfaces up to dimension $14$,
without any indication of the postulation of points in the fibers over multiple points, though.

We should point out that the existence of surfaces  with a $4$-dimension trisecant lines locus and, more generally,
of  embedded projective $n$-folds with fibers of unexpected postulation for a  generic projection to ${\P}^{n+1}$ is the source of great difficulties
in applying the general projection method to establish good Castelnuovo-Mumford regularity upper-bounds. For a discussion of the general projection
method, see Kwak's article \cite{Kw}.
This problem was raised by Greenberg \cite{Gre}, who showed how to bypass this issue for smooth surfaces, in view
of establishing Castelnuovo's regularity bound (Pinkham \cite{Pi}, Lazarsfeld \cite{L}).

A dimension count suggests that the trisecant lines locus of a surface in ${\P}^5$ should be $3$-dimensional;
so, surfaces in ${\P}^5$ with a $4$-dimensional trisecant lines  locus are said to have \textit{an excess of trisecant lines}.

Very few examples of surfaces with an excess of trisecant lines are  known and no classification has been established so far.
The first example of a surface with an excess of trisecant lines, known as
the \textit{Segre polygonal  surface} in the literature,
was constructed by  Segre \cite{S} in 1924.  It belongs to the family of White surfaces, which was constructed by
White and Gambier, also in 1924 \cite{G,W}. A modern reference is the paper of Gimigliano \cite{Gi}.

In his thesis \cite{D} , Dobler shows that the  Segre polygonal surface is the only polygonal surface with a $4$-dimensional trisecant lines locus.
Moreover, he shows  that the Segre polygonal surface is a degeneration of  smooth special Enriques surfaces in ${\P}^{5}$.
These surfaces are the only smooth surfaces in ${\P}^{5}$ with an excess of trisecant lines
(Conte and Verra \cite{CV}, Dolgatchev and Reider \cite{DR}).

Since polygonal surfaces belong to the family of White surfaces, it is natural to ask whether there exist other White surfaces
with an excess of trisecant lines.

The purpose of this work is  to show that the only White surfaces in ${\P}^{5}$ with an excess of trisecant lines are the Segre
polygonal surfaces.
To achieve this, we use a remarkable argument by Gambier \cite[p. 184-186 and p. 253-256]{G}.
Most of this work is intended to give modern rigor to his approach. Furthermore, this result allows us to deduce
information on the generic point  of the principal component of $W_{18}[5]$, the subvariety of the Hilbert scheme of
$18$-tuples of ${\P}^2$ special in degree $5$.  We prove the following results.
\begin{thm} Let $S_5$ be a \textit{possibly singular} White  surface in ${\P}^{5}$ and $p$ a generic point on $S_5$.

  \begin{itemize}

\item[(i)] Through the point $p$ there passes at least one $3$-secant line.

\item[(ii)] The surface $S_5$ has an excess of trisecant lines, if and only if $S_5$ is a Segre polygonal surface.

\item[(iii)] If $S_5$ doesn't have an excess of trisecant lines,  then  exactly $6$ trisecant lines, counted with multiplicity, pass through $p$.
\end{itemize}

\end{thm}
As a corollary of (i) we obtain
\begin{cor} The generic point of the principal component of $W_{18}[5]$ is
a smooth uniform $18$-tuple of points in the plane not lying on any curve of degree $4$.
\end{cor}
\section{Notations and basic facts}
\subsection{Linear systems of plane curves of given degree with assigned base points}

Let $S$ denote either ${\P}^{2}$ or a surface obtained from ${\P}^{2}$ by a finite succession of  blow-up at a single point, called $\sigma$-process.
Let $S\xrightarrow{\phi}{\P}^{2}$ denote the composition of these $\sigma$-processes.  If $S={\P}^{2}$, we set $\phi =id_{{\P}^{2}}$.
Let $L$ denote the linear equivalence class of a line in ${\P}^2$. A curve on $S$ is of degree $d$ if it is linearly
equivalent to $d\phi^{*}L$. Let us recall briefly the facts regarding linear systems of curves of degree $d$ with assigned base points and multiplicities
on $S$.

\begin{defn}[$n$-tuples of points on $S$] An $n$-tuple of points $P$ on $S$ is the data of
\begin{enumerate}
\item a sequence of distinct points on $S$, $\{p_1,\cdots ,p_{k}\}$ (These points
are the \underline{base points} of $P$.)

\item and a sequence of non-negative integers $\{ m_1,\cdots, m_k \}$ (These integers are the \underline{multiplicities} of $P$).

\end{enumerate}
such that $\sum_{i=1}^{k} m_{i}=n$.
\end{defn}

If $n=k$, we say that $P$ consists of distinct points. Let $P$ be a $n$-tuple of points on $S$.
Let $\tilde{S}\xrightarrow{\pi}S$, denote  the blow-up of $S$ at the base points of $P$. For $i\in \{ p_1,\cdots ,p_k \}$,
let $E_i$ denote the component of the exceptional locus of $\pi$ contracted to $p_i$ by $\pi$.
The curves on $S$ whose strict transform by $\pi$ are linearly equivalent to $d (\phi\circ \pi)^{*}(L)-\sum_{i=1}^{k}m_{i}E_i$,
form a linear system:\emph{the linear system of curves of degree $d$ on $S$ passing through $P$}. The plane curves of this system
are the vanishing locus of a polynomial in $H^{0}({\P}^2, {\I}_{P|{\P}^2}(d))$; thus, they are parameterized by points of the projective space
$|{\I}_{P}(d)|:={\P}(H^{0}({\P}^2, {\I}_{P|{\P}^2}(d)))$.
If the linear system $|d (\phi\circ \pi)^{*}(L)-\sum_{i=1}^{k}m_{i}E_i |$ is base-point free,
we say that the linear system $|{\I}_{P}(d)|$ is \emph{complete}.

Recall that the {\it virtual dimension } of $|{\I}_{P}(d)|$ is given by
\[
v({\I}_{P}(d))=\frac{d(d+3)}{2}  -\sum_{i=1}^{k} \frac{m_{i}(m_{i} +1)}{2}.
\]
We denote by $s$ the irregularity of the linear system $|{\I}_{P}(d)|$
\[
s:=dim H^{1}({\P}^2,{\I}_{P}(d))=dim (|{\I}_{P}(d)| )-Max(v({\I}_{P} (d)),-1).
\]
The linear system $|{\I}_{P}(d)|$ defines a rational  surface $X$ in ${\P}^{N}=|{\I}_{P}(d))|^{\vee}$,
image of $S$ by the rational map
\begin{equation}\label{Eq:1}
\begin{matrix}
\Phi \, : \, S \rightarrow {\P}^{N} \\
q\mapsto |{\I}_{P\cup q}(d)|^{\vee}
\end{matrix} 
\end{equation}
The surface $X$ is of degree ${\I}_{P}(d)^{2}$, the intersection number of two curves of the linear system $|{\I}_{P}(d)|$, 
and of sectional genus $\pi (X) =g_{{\I}_{P}(d)}$, the arithmetic genus of a  curve in $|{\I}_{P}(d)|$.

Following Dobler, we use a non-standard definition of a trisecant line; it is not defined as a line meeting $S$ along a $0$-dimensional 
scheme of degree at least $3$.
\begin{defn}[Trisecant lines]\label{def:tris} 
A line $l$ of ${\P}^N$ is a trisecant line to $X$ if $l$ does not lie on $X$ and cuts $X$ at (at least) three distinct points.
\end{defn}
\emph{The trisecant lines locus} of $X$ is then the closure in ${\P}^N$ of the union of all the trisecant lines to $X$.
 
A generic  trisecant line to $X$ passing through a generic point $\Phi(q)$ of $X$
can be directly seen on the linear system $|{\I}_{P}(d)|$. Indeed, such a line corresponds to a pair of distinct points
$\Pi:=(\pi_{1},\pi_{2})$ in $S$ such that the sublinear system of  $|{\I}_{P}(d)|$ of curves containing  $P\cup \Pi\cup \{q\}$
is $1$-irregular. We call such a pair $\Pi$ \emph{an associated pair to $|{\I}_{P}(d)|$ at $q$}.
Note that, with our special definition of trisecant lines, an associated pair to $|{\I}_{P}(d)|$ at $q$ defines a trisecant line to $X$
at $\Phi(q)$, if and only if the three points $\Phi(q),\Phi(\pi_1),\Phi (\pi_2)$ are distinct on $X$.

\subsection{Linear systems of plane curves of given degree through fixed intersection cycles}\label{sec:1}

Let $D_1$ and $D_2$ be two irreducible plane curves. Consider $Q:=D_1 \cap D_2$ the $0$-scheme of intersection of these two curves.
Let $\{ p_1,\cdots , p_m\}$ be the support of $Q$ . For every point  $p_i\in Supp Q$, let $n_i:=\mu_{p_i}(D_1,D_2)$ denote the multiplicity
of intersection of $D_1$ and $D_2$ at the point $p_i$. Then $Q$ gives rise to the $0$-cycle of the plane $\sum_{i=1}^{m} n_i p_i$, 
that we still denote by $Q$. We say that $Q$ consists of distinct points if $n_i=1$ for all $i=1,\cdots ,m$.
\emph{The linear system of plane curves of degree $d$ passing through the intersection cycle $Q$}
is ${\P}(H^{0}({\P}^2, {\I}_{Q|{\P}^2}(d)))$, where ${\I}_{Q|{\P}^2}$ is the ideal sheaf defining the $0$-scheme $Q$ in ${\P}^2$.
If $Q$ consists of distinct points, this linear system is just the linear system of curves of degree $d$ passing through the $p_i$'s.

Assume that $Q$ does not consist of distinct points; we may assume that $n_i\geq 2$ for $i=1,\cdots k$ and $n_i=1$ otherwise.
Suppose, for simplicity, that both curves $D_1$ and $D_2$ are smooth. 
A plane curve $C$ of degree $d$ belongs to $\mathcal{L}$ if and only if it contains the $0$-scheme $Q=D_1\cap D_2$.
Let $S^{(1)}\xrightarrow{\phi_1}{\P}^2$ denote the blow-up of ${\P}^2$ at $\{ p_1,\cdots ,p_k \}$. 
For $i=1,\cdots , k$, we define $E^{(1)}_i$ to be $\phi_1^{-1}(p_i)$. Denote by $Q^{(1)}$ (resp. $D_1^{(1)}$ and $D_{2}^{(1)}$) the strict transform
of the scheme $Q$ (resp. $D_1$ and $D_2$) by $\phi_1$; then, $Q^{(1)}$ is the $0$-scheme of intersection on $S^{(1)}$
of the curves  $D_1^{(1)}$ and $D_{2}^{(1)}$. We have, moreover 
\[
\phi_{1}^{-1}(D_1)\cap \phi_{1}^{-1}(D_2)=Q^{(1)}\cup E_{1}^{(1)}\cup\cdots \cup E_{k}^{(1)}
\]
Thus, the strict transform $\mathcal{L}^{(1)}$ by $\phi_1$ of the linear system $\mathcal{L}$
is the system of curves  of $S^{1}$ linearly equivalent to 
\[
\phi_1^{*}({\O}_{{\P}^2}(d))-\sum_{i=1}^{k}E^{(1)}_{i}
\]
and containing the $0$-scheme $Q^{(1)}$.

Recall that the multiplicity at a point $p$ of an irreducible curve $C$ lying on a surface 
is the sum of the multiplicities of infinitely near points of $C$ of order $1$ over $p$ (see \cite[p.33]{Bea}).
Moreover, there is a point $q\in Supp(Q^{(1)})$ such that $\phi_1 (q)=p_i$ if and only if $n_i\geq 1$. 
Since we have assumed that both $D_1$ and $D_2$ are smooth, such a point $q$ is unique, if it exists; we denote it by $p_{i}^{(1)}$. 
The point $p^{(1)}_{i}$ is an infinitely near point of order $1$ over $p_i$ for both $D_1$ and $D_2$. 
Let $n_i^{(1)}:=\mu_{p_i^{(1)}}(D_{1}^{(1)} , D_{2}^{(1)}$.
If $Q^{(1)}$ does not consist of distinct points, repeat the process of blowing-up the multiple points of the support of $Q^{(1)}$.
This process stops after a finite number of steps (see \cite[theorem 4.2.5 ]{Na}). We get a sequence
of blow-up maps
\[
S^{(l)}\xrightarrow{\phi_l}S^{(l-1)}\xrightarrow{\phi_{l-1}} \cdots S^{(1)}\xrightarrow{\phi_1} {\P}^2,
\]
such that the intersection cycle $Q^{(l)}$ consists of distinct points and
the strict transform by $\phi:=\phi_{1}\circ\cdots\circ\phi_{l}$ of the linear system $\mathcal{L}$ is the linear system of curves in 
\[
|\phi^{*}({\O}_{{\P}^2}(d))-\sum_{i=1}^{k} \sum_{j=1}^{n_{i}-1}E^{(j)}_i|
\]
passing through the distinct points of $Q^{(l)}$.
Indeed, there are only $n_i -1$ cycles among $Q^{(1)},\cdots ,Q^{(l)}$ whose support contains a point over $p$, for we have
\cite[p.33]{Bea}
\[
\mu_{p_i}(D_1,D_2)=\mu_{p_i}(D_1)\mu_{p_i}(D_2)+\sum_{x\in A(p_i)} \mu_x(D_1)\mu_x(D_2),
\] 
where $A(p_i )$ is the set of infinitely near points of both $D_1$ and $D_2$ over $p_i$. 
We can formally replace the cycle $Q$ by the cycle of ``distinct or infinitely near points''
$\sum_{i=1}^{k}(p_i+\sum_{j=1}^{n_{i}-1} p_i^{(j)})+\sum_{i=k+1}^{m}p_i$. 
The linear system $\mathcal{L}$ is thus the linear system 
of curves of degree $d$ passing through the proper distinct points $p_1,\cdots p_m$ and the infinitely near points
$p_1^{(1)},\cdots ,p_{1}^{(n_{1}-1)},\cdots ,p_{k}^{(1)},\cdots, p_{k}^{(n_{k}-1)}$. 
Of course, a similar analysis can be made to understand better  the linear system $\mathcal{L}$ in case $D_1$ or $D_2$ is singular. 
 
\subsection{Standard irregularity estimation techniques}
Let $S={\P}^2$. Recall the standard geometric interpretation of the irregularity of pencils.
Let $|{\I}_{P}(d)|$ be as in section \ref{sec:1}. Suppose that $|{\I}_{P}(d)|$ contains a pencil, \ie 
$dim H^0({\I}_{P}(d))\geq 2$. Suppose, moreover, that the linear system $|{\I}_{P}(d)|$ is complete and that all the multiplicities of $P$ are equal to one.

Any two distinct curves $D$ and $D^{\prime}$ of the linear system $|{\I}_{P}(d)|$ meet in a $0$-dimensional scheme.
Let $Q$ denote the cycle of intersection $D\cdot D^{\prime}\setminus P$. As $D^{\prime}$ varies in $|{\I}_{P}(d)|$, the cycles $Q$
fit into a linear series $\chi (D, |{\I}_{P}(d)|)$ on $D$, called \emph{the characteristic series} cut on $D$.
A member of this series is said to be \emph{residual to $P$} on $D$.

\begin{prop}[Duality Theorem ]\label{prop:duality}
Let  $|{\I}_{P}(d)|$ be as above. Let $D$  be a generic curve in $|{\I}_{P}(d)| $ and   $Y\in \chi( D, | {\I}_{P}(d)|)$. 
Denote by $s$ the irregularity of $|{\I}_{P}(d)| $. 
Then, $s=dim H^{0}({\P}^2, {\I}_{Y|{\P}^2 }(d-3))$, where  ${\I}_{Y|{\P}^2 }$ is the ideal sheaf defining $Y$ in ${\P}^2$.
\end{prop}

For a proof of this proposition see Griffiths and Harris' book \cite[p.713-716]{GH}.

In order to estimate the irregularity of $|{\I}_{P}(d)|$, we need  to control the configuration of the group of points residual to $P$
with respect to the linear system.
This is taken care of by the classical residuation theorem, which is a direct rewriting of Noether "AF+BG" theorem.
We require the stronger version of it, which allows singular points in the intersection cycles.
It can be found, for instance, in Walker's book on algebraic curves \cite[theorem 7.2]{Wa}.

\begin{prop}[residuation theorem]\label{prop:residu}
Let $C_n$ and $C_m$ be  plane curves of degree $n$ and $m$. Using cycle notation, we write $C_n \cdot C_m =P+Q$.
Suppose that there are two integers $n_{1}$ and $n_{2}$ such that $n_{1}+n_{2}\ge n$
 and there exist two curves of degree $n_1$ and $n_2$, respectively, such that:
\[
C_{n_1}\cdot C_{m}=P+P^{\prime}\,\,\text{and}\,\, C_{n_2}\cdot C_m =Q+Q^{\prime}.
\]

Moreover assume that
\begin{itemize}
\item[(1)] either $C_m$ is \emph{smooth} at any point of the support of $P+Q+Q^{\prime}+P^{\prime}$,

\item[(2)]or, if $n_2=0$ (classical residuation theorem), that $C_m$ is smooth at the points of $C_m \cdot C_n$.

\end{itemize}

Then there exists a curve of degree $n_1+n_2-n$ such that
\[
C_m \cdot C_{n_1+n_2 -n}=P^{\prime}+Q^{\prime}.
\]
\end{prop}

\subsection{White and polygonal surfaces}

Let $d\geq 5$ be an integer.

\begin{defn}[White surfaces]
Pick $\binom{d+1}{2}$ \emph{distinct} points $P$ in the plane, not lying on any curve of degree $d-1$;
so the linear system $|{\I}_{P}(d)|$ is complete and regular. The non-degenerate rational surface $S_d$ it defines in ${\P}^{d}$,
is called a \textit{White surface} (Gimigliano \cite{Gi}).
\end{defn}
We should stress that Gimigliano studies much more general White surfaces than the classical ones we are interested in.

\begin{defn}[polygonal surfaces]
A simple way to construct $\binom{d+1}{2}$ points of ${\P} ^{2}$ not lying
on any curve of degree $(d-1)$ is to
take the points of intersection of $(d+1)$ lines general enough.
White surfaces obtained this way are called \textit{polygonal surfaces}.
\end{defn}

\begin{defn}[Segre polygonal surfaces]
A Segre polygonal surface is a polygonal surface obtained by imposing that the $(d+1)$ lines defining the polygonal surface are tangent to a fixed conic.
\end{defn}

We refer the reader to Gimigliano's paper \cite{Gi} for general properties
of White surfaces. We gather here the properties we need further on.

\begin{prop}[Gimigliano \cite{Gi}]Let $S_d$ be a White surface and $P$ the points in ${\P}^2$ used to construct $S_d$.
Let $\tilde{\P}$ be the blow-up of the plane along $P$ and $\phi$ the map $
{\tilde{\P}}\xrightarrow{} {\P}^{d}$
defined by $|{\I}_{P}(d)|$.
\begin{itemize}
\item[(i)] If no $d$ points of $P$ are aligned, the map $\phi$ is an embedding.

\item[(ii)] Under the map $\phi$, a line $l$ containing $d$ of the points $P$ is contracted to a $(d-1)$-fold point.
Any two such lines have to meet in a point of $P$.

\item[(iii)] All the singularities of $S_d$ are obtained by contraction of
lines containing $d$ base points.

\end{itemize}
\end{prop}

For a proof, see Gimigliano \cite{Gi} and Dobler \cite{D}.

Dobler deduces the following result on trisecant lines to \emph{smooth}
White surfaces from Bauer's results on inner projections \cite{Ba}.

\begin{prop}[Dobler, \cite{D} proposition 2.7 ] \label{prop:Dob1} Let $S_5$
   be a smooth White surface in ${\P}^5$ and $q$ a generic point of $S_5$.
Then, there exists a trisecant line to $S_5$ passing through $q$.
\end{prop}

Dobler also proves the following results concerning the trisecant locus of
polygonal surfaces.

\begin{prop}[Dobler, \cite{D}  proposition 3.17] \label{prop:Dob2} Among polygonal surfaces in ${\P}^d$, only the Segre polygonal surfaces
have a $4$-dimensional trisecant locus. In fact, for any non-Segre polygonal surface, there are at most finitely many trisecant lines to $S_d$.
For a general polygonal surface there are none.
\end{prop}

\subsection{The variety parameterizing $N$-tuples of points special in degree $d$}

In this section we briefly sketch the results regarding the (generally non-irreducible) varieties $W_{N} [d]$ that we will use later on;
for further details, we refer the reader to Coppo's paper \cite{C} and Ellia and Peskine's paper \cite{EP}.

Following Ellia and Peskine, we call group of points a $0$-dimensional scheme of ${\P}^2$. A group of points $Z$
is \emph{special in degree d} if $H^{1}({\P}^2, \I_{Z|{\P}^2}(d))\not= 0$.
The groups of points of degree $N$ (called $N$-tuples of points in Coppo's terminology) special in degree $d$ 
are parameterized by a possibly non-irreducible subvariety, $W_N [d]$, of the Hilbert scheme $Hilb^N({\P} ^2)$ of groups of points of length $N$.

Let $Z$ be such a group of points.  Choose a \textit{generic} line $l$ in ${\P}^{2}$, so that it avoids the support of $Z$.
We may assume its  equation is given by $z=0$, where $x,y,z$ are  projective  coordinates of the plane.
Let us recall that the coordinate ring $A_Z$ of  a $0$-dimensional scheme $Z$ in ${\P}^2$
has a graded $A=\hbox{\pr C}[x,y]$-module structure, which has a $A$-module resolution of the form
\[
0\rightarrow\oplus_{i=0}^{s-1}A[-n_i]\rightarrow \oplus_{i=0}^{s-1}A[-i
] \rightarrow A_Z\rightarrow 0
\]
with $ n_0\geq n_1\cdots\geq n_{s-1}\geq s$.
These positive integers $n_i$ give a numerical invariant for $Z$ introduced by Gruson and  Peskine \cite{Pes},
the \emph{numerical  character} $\chi  (Z)= (n_0,\cdots ,n_{s- 1})$.
The integer $s$, called the \emph{height} of $\chi (Z)$, is the minimal degree of a plane curve containing $Z$.
The  \emph{index of specialty} of $Z$, $n_0 -2$, is the greatest integer $n$ such that $H^{1}({\I}_{Z}(n)) \not= 0$.
For the moment, assume that the linear system $|{\I}_Z(d)|$ is complete.

\emph{If $Z$ is generic} among the $N$-tuples of points of character $\chi (Z)$,
the superabundance of $|{\I}_Z (d)|$ can be computed  directly from $\chi (Z)$ by the following formula \cite[\S 1.1]{C}
\[
h^1(\I_Z (d))=\sum_{i=0}^{s-1} (n_i-d-1)^{+}-(i-d-1)^{+},
\]
where $(n)^{+}$ is  zero, if the integer $n$ is negative, and $n$ otherwise.

Of course, the degree of $Z$ can be recovered from $\chi (Z) $:
\[
deg(Z)=deg(\chi ):=(\sum_{i=0}^{s-1} n_i)-s-1.
\]
The following proposition is useful to compute $\chi$ from the geometry of
the group of points $Z$.

\begin{prop}[Davis \cite{Da}, Ellia-Peskine \cite{EP}]\label{prop:EP} Let
$\chi =(n_0,\cdots ,n_{s-1})$ be the
numerical character of a group of points $Z$. Suppose that for some index $
1\leq t\leq s-1$
we have $n_ {t-1}>n_t +1$. Then there exists a curve $T$ of degree $t$ such that,
if
\begin{enumerate}
\item $Z^ {\prime}$ is the group of points $Z\cap T$,

\item $f=0$ is the local equation of $T$ in $A$ and

\item $Z^{\prime\prime}$
is the residual group of points of $Z^{\prime}$ in $Z$ with respect to $T$,

\end{enumerate}
we have $A_{Z^{\prime\prime}}=fA$ and $\chi (Z^ {\prime\prime})=(m_i)$
with $m_i=n_{t+i}-t$.
\end{prop}

\section{Projection of a surface from a multisecant line}

Let  $H$ denote the class of hyperplanes in ${\P}^5$. Classically, we get geometric information on a surface $X$ 
from the study of its double locus $D$ by a generic projection to  ${\P}^{3}$.
If $X$ is normal of degree $d$ in ${\P}^{r}$ and $K_{X}$ is its canonical divisor, then a well-known
consequence of the theory of subadjoint systems is that the class of $D$ in $Pic(X)$ is $(d-4)H-K_{X}$ (\cite{R}, \cite{Z}).

We need to establish a similar formula for the double locus of the projection of a smooth surface $X$,
if the center of projection is a multisecant line.
This can be done by viewing  this  projection as a limit of regular projections, following a method  
that Franchetta (\cite{Fr}, \cite{Ci}) used to prove  his famous theorem on the irreducibility of the double locus
 of a generic  projection of a codimension two surface.

\begin{prop}\label{prop:proj} Let $X$ be a smooth, non-degenerate, surface
of degree $d$ in  ${\P}^5$. Assume, moreover, that $l$ is a line in ${\P}^5$ cutting $X$ along a 
$0$-dimensional scheme of multiplicity  $\delta$.
Let us consider the projection ${\pi}_{0}$ of $X$ to  ${\P}^{3}$ with center of projection $l$.
Suppose that  any plane containing $l$ intersects $X$ at most in a finite number of points.

Then the class, in $Pic(X)$, of the double locus $D$ of $\pi$ is given
by
\begin{equation}
D=(d-4-\delta )H-K_X . \label{Eq:double locus}
\end{equation}
\end{prop}

We believe it is a classical result; although, we could not find any reference.

\begin{pf} If the plane $<q,l>$ meets $X$ for all points $q$ in ${\P}^{5}\setminus l$, then ${\pi}_{0}(X)={\P}^{3}$.
Therefore, we can find a plane  $\Pi$ containing $l$ such that $\Pi\cap X=l\cap X$.
First, let us construct a family of regular projections degenerating at ${\pi}_{0}$.

Consider the pencil $\mathcal{P}$ of lines  contained in $\Pi$ and passing
through a given generic point $p_0$ on $l$.
Let $0$ be a generic point on $l$. Let $\Delta$ be a generic line meeting $l$ at $0$; then,
$\Delta$ parameterizes the lines $l_t$ of the pencil
${\mathcal{P}}$. Let us denote by $\{ p_{1},\ldots ,p_{s}\} $ the support
of $l\cap X$ and by $m_{i}$ the multiplicity of intersection of $l$ and $X$ at the point $p_{i}$.

For $t\not =0$, the projection centers $l_t$ do not meet $X$, so the projections $\pi_{t}$ of ${\P} ^{5}\times \{t\}$ onto ${\P}^{3}\times \{ t\}$
define a rational map  $\pi$ from the $3$-fold $X\times \Delta$ to $\mathcal{F}=\pi (X\times\Delta)$  in ${\P}^{3}\times \Delta$,
whose indeterminacy locus is $\{ (p_1,0),\ldots ,(p_{s},0)\} $.

Thus, the blow-up of $X\times \Delta$ at the points $(p_1,0),\cdots ,(p_{s},0)$, $\phi\, :\, \mathcal{ X}\rightarrow X\times \Delta$,
induces a regular map $g$ fibered over $\Delta$, which resolves the indeterminacy of $\pi$.

We have the following diagram:
\[
\xymatrix{
\mathcal{X}\ar[d]^{\phi}\ar[dr]^{g}& \\ {X\times \Delta}\ar@{-->}[r]^{\pi}\ar[d]^{
h} & {\mathcal{F}\subset {\P}^{3}\times \Delta}\ar[dl]^{\zeta}\\ \Delta &
}
\]
where $h$ (resp. $\zeta$) denotes the projection to the second factor.
For $1\leq i\leq s$, let $E_i$ denote the exceptional divisor of $\mathcal X$ over $(p_i,0)$.

Since $h\circ \phi $ is a morphism from  $\mathcal{X}$ onto a $1$-dimensional smooth base $\Delta$, it is a flat morphism \cite[III 9.7 p.257]{H}.

Its special fiber $(h\circ\phi)^{-1}(0)$ is simply given by
\[
{\mathcal X}_{0}={\tilde X}\times \{ 0\} \cap_{i=1}^{s} E_i
\]
where $\tilde X $ denotes the blow-up of $X$ at $p_ {1},\ldots , p_{s}$.
Recall that the double locus $D$ of the map $g$ is defined  on $\mathcal{F}
$ by the ideal sheaf of the conductor
$\hbox{\got c}:=({\O}_{\mathcal {F}} \colon g_{*}({\O}_{\mathcal{ X}}))$.

Since $g$ is a finite morphism, $\hbox{\got c}$ is also an ideal  sheaf over ${\O}_{\mathcal{X}}$, defining the double locus
$D^{\prime}=g^{-1}(D)$ of $g$ on $\mathcal{X}$, see for instance \cite{R}.
The second projection from $\mathcal{ F}$ onto  $\Delta$ is a flat morphism, so
  the   ${\O}_{{\mathcal F}_{t}}$-sheaf ${\O}_{\Delta,t}$ is  locally free.
If we tensor the previous exact sequence by this sheaf, we simply get the exact sequence
defining the double locus $D_{t}$ of the map $g_{t}$,
restriction of $g$  to the fiber $\mathcal{X}_{t}$. Thus, in each  fiber $\mathcal{X}_{t}$,
the double locus $D$ of $g$ restricts to the double locus of $g_t$.
Since $\hbox{\got c}$ is an invertible ${\O}_{\mathcal X}$-module, $(D^{\prime}_{t})$
forms a flat family of divisors over $\Delta$.
For $t\not= 0$, $\mathcal{X}_{t}=X\times \{ t\}$ and  $g_t=\pi_{t}$,
the class of $D^{\prime}_{t}$ in $Pic(\mathcal{X}_{t})$ is  classically given by \cite{R, Z}
\[
D^{\prime}_{t}=(d-4)H^{*}_{t}-K_{\mathcal{X}_{t}}
\]
where $H^{*}_{t}=g_{t}^{*}((\zeta^{*}(H))|_{\mathcal{F}_{t}}) $, $H$ being the class of hyperplanes in ${\P}^{3}$.
The divisor $H^{*}_{t}$ is naturally the restriction of $(g)^{*}({\zeta}^{*
}(H)|_{\mathcal{X}})$ in $Pic(\mathcal{X})$.
\begin{claim}
For $t\not= 0$, we have $K_{\mathcal{X}_{t}}=K_{\mathcal{X}}|_{\mathcal
{X}_{t}}$.
\end{claim}
\begin{pf}
Since $\mathcal{X}$ is a blow-up of $X\times \Delta$, we have $K_{\mathcal{X}}=\phi^{*}(K_{X\times \Delta})+\sum_{i=1}^{s} E_{i}$.
Since, for all $i=1,\ldots ,s$,  we have $\mathcal{X}_{t}\cap E_{i}=\emptyset$,
we find  $K_{\mathcal{X}}|_{\mathcal{X}_{t}}
=\phi_{t}(K_{X\times \Delta}|_{X\times\{ t\} })$. Recall that
\[
Pic(X\times \Delta)\simeq \eta^{*}(Pic(X))\oplus h^{*}(Pic(\Delta))\simeq \eta^{*}(Pic(X))\oplus h^{*}(H_{\Delta}{\Z})
\]
where $\eta$ is the projection from $X\times \Delta$ onto  the first factor
  and $H_{\Delta}$ stands for the class
of points in $Pic(\Delta)$. From $K_{X\times \Delta}=\eta^{*}(K_{X}) +h^{*}(-2 H_{\Delta})$ we deduce
that $K_{X\times\Delta}|_{X\times\{ t\}}=\eta^{*}(K_{X})|_{X\times \{ t\}}$.
\qed
\end{pf}
By \cite[proposition 9.7 p.258]{H}, the family $(D^{\prime}_ {t})_{t\not=_0}$
has a unique flat limit, so the double locus $D^ {\prime}$ of $g$ on
$\mathcal{X}$  has class $(d-4)H^{*}-K_{\mathcal{X}}$. Finally, we must determine  the double locus  divisor of $g_{0}|_{\tilde{X_{0}}}$.
Note that $X\times  \Delta$ is embedded  in ${\P}^{5}\times \Delta$ by the
linear system $H^{\prime} =h^{*}(H_{X})$;
so the linear system defining  $g$, is $|H^ {\prime}-\sum_{i=1}^{s} m_{i}  E_{i}|$ and $g_0$ is defined by the restriction
of this  system to $\tilde{  X_{0}}$.
Since  $K_{{\mathcal X}_{0}}=\phi_{0}^{*}(K_X\times \{ 0\} )+\sum_{i=1}^{s}E_i$,
we get $D_0^{\prime}=(d-4)g_{0}^{*}(H_{0})-\phi_{0}^{*}(K_X\times \{  0\} )-\sum_{i=1}^{s}E_i$.
Moreover, the restriction of  $|H^{\prime}- \sum m_{i}  E_{i}|$ to $X\times \{0\}$ embeds $E_{i}$ as a plane in ${\P}^{3}\times \{ 0\} $.
The curves $e_i=\tilde X\cap E_i$ for $1\leq i\leq s$  sit in the  double locus of $g_0$ and are mapped by $g_0$
to  hyperplane sections of $\mathcal{F}_0$.
Thus, the  restriction  of the double locus of $g_0$ to ${\tilde X} _{0}$ can only differ from the double locus of $g_{0}|_{\tilde X}$
along the curves $e_i$. Besides, the double locus of $g_{0}|_{\tilde X}$ is not supported along these curves.
Therefore, the class of the double locus of the projection of $X$ from $l$
is given by the formula
\[
\phi_{0*}((d-4)g_{0}^{*}(H_0)-\phi_{0}^{*}(K_{X\times \{ 0\}}) -\sum_{i=1}^{s} e_i-g_{0}^{*}(\delta H)).
\]
  \qed
\end{pf}

\begin{rem} Let $S$ be a singular surface of ${\P}^{5}$, obtained as the image of ${\P}^{2}$ by a complete linear system of plane curves
with assigned base points. If $S$ satisfies all the hypothesis of the theorem but smoothness, a slight modification of the argument shows
that the same formula holds for the double locus in the Picard group of $\tilde{ \P}$, the blow-up of ${\P}^{2}$ at the base points.
\end{rem}

\section{Existence of a trisecant line through a  generic point of $S_5$}

Following Gambier \cite[p. 184-186]{G}, we extend Dobler's result on the existence of trisecant line to the case of singular White surfaces. 
We present Gambier's beautiful construction in modern language and show how to fill a few gaps.

\begin{thm}[Gambier]\label{thm:1} Let $P$ be a finite set of $15$ distinct
points in ${\P}^{2}$, such that $P$ is not contained in any curve of degree $4$. Pick a generic point $q$ in ${\P}^{2}$.
Then, there exists an associated pair $\Pi$ to $|{\I}_{P}(5)|$ at $q$.
\end{thm}
\begin{pf}  Let $D$ be  a generic  curve of $|{\I}_{P+q}(5)|$, so that $D$ is smooth. Fix a  pencil
$<D,D^{\prime}>$ of curves in $|{\I}_{P+q}(5) |$.
Let $Q$ denote the $0$-scheme of length $9$ residual to $P+q$ in the intersection cycle $D\cdot D^{\prime}$.
Pick a basis $D,D_1,\cdots , D_4$ of the linear system $|{\I}_{P+q}(5)|$ and choose $D^{\prime}$ generic in
the linear subsystem $<D_1,\cdots ,D_4>$, so that $D^{\prime}$ is smooth.

An easy dimension count shows that the virtual dimension
of the linear system of cubic curves passing through the cycle $Q$ is non negative.

\begin{lem}\label{lem:unique}
Let $D$ and $D^{\prime}$ be as above, then there is a unique cubic curve passing through $Q$.
\end{lem}
\begin{pf}
Suppose, to the contrary, that there is a pencil of cubic curves, $<C,C^{\prime}>$,  passing through $Q$.
We have only two cases to consider:
\begin{itemize}
\item[(A)] either those cubic curves are generically irreducible,

\item[(B)] or there is a fixed line or a fixed conic in the pencil $<C,C^{\prime}>$.
\end{itemize}
Let us  prove the lemma in case (A).
First, assume that both cubic curves $C$ and $C^{\prime}$ are irreducible.
We claim that these curves are smooth along the supporting points of $Q$.

Let $r$ be a point of support of the cycle $Q$. We have, by construction,
\[
\mu_{r}(C,C^{\prime})\geq \mu_{r}(D,D^{\prime}),
\]
where $\mu_r (C,C^{\prime})$ denotes the multiplicity of intersection of $C$ and $C^{\prime}$ at $r$.
Since the two cubic curves are assumed irreducible, by summation of the previous inequalities,
we get
 \[
9\geq \sum_{r\in Supp Q}\mu_{r}(C,C^{\prime})\geq \sum_{r\in Supp Q} \mu_{r}(D,D^{\prime})= 9,
\]
so for every $r\in Supp Q$, we have
\[
\mu_{r}(C,C^{\prime})= \mu_{r}(D,D^{\prime})\geq \mu_{r}(C)\mu_{r}C^{\prime}.
\]

If $r$ is not a multiple point of $Q$, the curves $C$ and $C^{\prime}$ are therefore smooth at the point $r$.

Suppose, now, that $r$ is a multiple point of $Q$. We have $\mu_{r}(C,C^{\prime})= \mu_{r}(D,D^{\prime})$,
so we get
\begin{equation}\label{Eq:e1}
\mu_{r} (C)\mu_{r} (C^{\prime})+\sum_{x \in A(r)} \mu_{x}(C)\mu_{x}(C^{\prime})=\mu_{r} (D)\mu_{r} (D^{\prime})+\sum_{x \in B(r)} \mu_{x}(D)\mu_{x}(D^{\prime}),
\end{equation}
where $A(r)$ (resp. $B(r)$) is the set of infinitely near points of both $C$ and $C^{\prime}$ (resp. $D$ and $D^{\prime}$ ) over $r$.
From the analysis made in section \ref{sec:1}, we find $B(r)\subset A(r)$. Since $D$ and $D^{\prime}$ are smooth, we find $\mu_{r}(D)=\mu_{r}(D^{\prime})=1$
and $\mu_{x}(D)=\mu_{x}(D^{\prime})=1$  for all $x\in B(r)$. Thus, from equation \ref{Eq:e1}, we deduce that $A(r)=B(r)$, $\mu_{r}(C)=\mu_{r}(C^{\prime})=1$ 
and $\mu_{x}(C)=\mu_{x}(C^{\prime})=1$  for all $x\in A(r)$. In particular, $C$ and $C^{\prime}$ are smooth at $r$.

Since $C$ and $C^{\prime}$ are smooth at the supporting points of $Q$, we may apply the residuation theorem for $n=m=3$, $n_{1}=5$
and $n_{2}=0$. We find  that the cycle of length $6$, $V$, residual to $Q$ on $D\cdot C$, is
the complete intersection of $C$ with some conic $C_2$.
Applying, now, the residuation theorem to $V$ and $Q$ with $n_{1}=5$, $m=5$, $n=3$ and $n_{2}=2$,
we deduce that $P+q$ lies on a quartic curve. This contradicts our assumption on $P$.

Only remains to prove the lemma in case (B). Either the pencil of cubics is
the composite of a line $l$ with a pencil
of conics or a fixed conic $C_0$ composed with a pencil of lines.
Assume that we are in the first case and that the generic conic is irreducible, hence smooth. Thus, we can
pick two smooth generators of the pencil: $C_2$ and $C_2^{\prime}$. Applying the residuation
theorem as before with $C_2\cup l$ as cubic curve and $C_2^{\prime}$ as conic, we get the same contradiction.
The case of a fixed conic and a pencil of lines is similar.
\qed
\end{pf}

The following proposition  and lemma are a key ingredient used by Gambier without a proof.
\begin{prop}\label{prop:19}
The cubic curves $C$,  constructed this way, vary in a linear system as the
quintics $D^{\prime}$ vary in $|{\I}_{P+q}(5)|$.
\end{prop}

\begin{pf}  
It is enough to show that for any pencil $\mathcal{P}$ of quintic curves in $|{\I}_{P+q}(5)|$, such that $D\not\in\mathcal{P}$,
the family of cubic curves passing through the residual groups of points to $|{\I}_{P+q}(5)|$ on $D$ vary in a rational pencil.
\cite[p. 25]{Z}
As $t$ varies in ${\P}^1$, the parameter space of the pencil $\mathcal{P}$, the curves $\mathcal{D}_{t}$ of the pencil vary in a rational
family $\mathcal{D}$ over ${\P}^1$, sub-family of the trivial family ${\P}^2\times{\P}^1\xrightarrow{g}{\P}^1$.
We have 

 \[
\xymatrix{ 
 \mathcal{D} \subset & {{\P}^2 \times {\P}^1} \ar[d]^{f} \ar[dr]^{g} & \\
& {\P}^2 & {\P}^1 \\
}
\]
where $f$ (resp. $g$) is the projection onto the first factor (resp. the second factor) of ${\P}^2 \times {\P}^1$; 
the map $g$ is flat, since $g$ is surjective  and ${\P}^1$ is a smooth $1$-dimensional variety.
Consider the family 
\[
\mathcal{Q}:=\mathcal{D}\cap (D\times {\P}^1 )\setminus (\{P+q \}\times {\P}^1)\xrightarrow{}{\P}^1.
\]
For $t$ generic in ${\P}^1$, both $\mathcal{D}_{t}$ and $D$ are smooth; thus, by the previous lemma, there is a unique cubic curve $\mathcal{C}_t$
passing through the intersection cycle $\mathcal{Q}_{t}$.
A bigraded free resolution of ${\O}_{\mathcal Q}$ as ${\O}_{{\P}^2\times {\P}^1}$-module is of the form 
\[
\cdots \xrightarrow{}\oplus_{i=1}^{k}f^{*}({\O}_{{\P}^2}(-n_i))\otimes g^{*}({\O}_{{\P}^1}(-m_i))\xrightarrow{}{\O}_{{\P}^2\times{\P}^1}\xrightarrow{}{\O}_{\mathcal Q}\xrightarrow{} 0,
\]
for some non-negative integers $k$, $m_i$ and $n_i$.
Tensoring this resolution by ${\O}_{{\P}^2\times \{t\}}$, we get a presentation of ${\I}_{{\mathcal Q}_{t}}$
\[
\xrightarrow{}\oplus_{i=1}^{k}{\O}_{{\P}^2}(-n_i)\xrightarrow{}{\I}_{{\mathcal Q}_t}\xrightarrow{} 0.
\]
By lemma \ref{lem:unique}, (for $t$ generic) there is an index $i_{0}$ in $\{ 1,\cdots ,k\}$ such that $n_{i_0}=3$.
Thus, there exist a non-negative integer $n$ and a non-zero polynomial $\sigma$ of bidegree $(3,n)$ whose vanishing locus contains the scheme $\mathcal{Q}$. 
For $t$ generic in ${\P}^1$, the curve $(\sigma (t)=0)$ is the unique cubic curve passing through the $0$-scheme $\mathcal{Q}_t$.
Since the plane curves of the family $(\sigma=0) \xrightarrow{g}{\P}^1$ have the same degree,
hence the same Hilbert polynomial, $g$ is flat. Thus, $(\sigma =0)\xrightarrow{g}{\P}^1$ is a rational pencil.

\qed
\end{pf}

Let us choose a basis $D,D_1,\cdots ,D_4$ of $|{\I}_{P+q}(5)|\simeq {\P}^4$. As a consequence of the previous proposition, we get:

\begin{lem}\label{lem:residu}
The characteristic series  of  $|{\I}_{P+q}(5)|$  defines an algebraic injective morphism
$\Psi: {\P}^{3}=<D_1,\cdots ,D_4>\xrightarrow{} {\P}^{9}=|{\O}_{{\p}^{2}}(3)|$ which  sends a curve $C$ in $<D_1,\cdots ,D_4> $
to the \emph{unique} curve  of degree $3$ passing through the cycle $(D\cdot C)-(P+q)$.
\end{lem}

\begin{pf}
Only remains to prove the injectivity of $\Psi$. Suppose, to the contrary,
that there  exist two quintics $C$ and $C^{\prime}$ of the linear system $<D_1,\cdots ,D_4>$, defining the same cubic $C_3$.
We write
\[
C\cdot D=P+q+R\,\,\,\text{and}\,\, C^{\prime}\cdot D=P+q+R^{\prime}
\]
\[
D\cdot C_3=R+V=R^{\prime}+V^{\prime}.
\]
It follows that $R$ and $R^{\prime}$ contain at least $3$ common points.

Assume that $R$ and $R^{\prime}$ contain exactly
$3+n$ common points $\Pi$, for $0\leq n\leq 9$. The linear system $|\I_{P+q+\Pi}|$ is therefore special.
Let us write $R=\Pi +W$ and $R^{\prime}=\Pi +W^{\prime}$; we have $V=W^{\prime}+S$ and $V^{\prime}=W+S$, 
where $S$ is a $(n-3)$-tuple of points.
Pick $\Pi_0$ in $\Pi$ such that ${\Pi}_0 +W$ has length $6$.
By the duality theorem, $W$ (resp. $W^{\prime}$) lies on some conic $C_2$ (resp. $C_2^{\prime}$).
Finally, we apply the residuation theorem with $n=5$, $m=3$, $n_1=5$
and $n_2=2$,
to the curves $D$, $C_3$, $C$ and $C_2^{\prime}$.
Thus, $W+W^{\prime}+2\Pi_0$ lies on some conic $C^{\prime\prime}_2$.
Thus, we have $C^{\prime\prime}_{2} \cdot D\supset W+W^{\prime}+2\Pi_0$;
this gives a contradiction,
since $D$ is irreducible and $W+W^{\prime}+2\Pi_0$ has length $12$.
\qed
\end{pf}

We are now able to prove the theorem. The locus in ${\P}^{9}=|{\O}_{{\p}^{2}}(3))|$ of this family of cubic curves
is  the image of ${\P}^{2}\times {\P}^{2}\times {\P}^{2}$  by the Segre embedding. This is a variety of codimension $3$ in ${\P}^{9}$.
Since $\Psi$ is injective, we may choose $D^{\prime}$ such  that the unique cubic curve passing through $Q$,
the group of points of $D$ residual to $P+q$ in $D\cdot D^{\prime}$, is the union of $3$ lines $l_{1}$,  $l_{2}$ and $l_3$.
Let $T$ denote the triangle $l_{1}\cup l_{2}\cup l_3$.
Let $V^{t}$ denote the group of points on $D_t$ residual to $Q$ in $D_{t}\cdot T$, as $D_{t}$ varies in the pencil $<D,D^{\prime}>$.
We can find $t_{0}\in {\P} ^{1}$ such that $3$ points of $V^{t_{0}}$ lie on a line, say $l_{1}$.

Since $D$ is irreducible, there are no more than $3$ points of $V_{t_{0}}$, counted with multiplicity on $l_1$.
Indeed otherwise, there would not be a unique cubic curve passing through $Q$.
As $t$ tends to $t_{0}$,  $D_{t}\cdot l_1$ tends to a length $5$ cycle, $\Pi +(V_{t_{0}}\cdot l_1 )$.
Since $\Pi$ is a sub-cycle of $Q$, $\Pi$ belongs to the intersection cycle $D\cdot l_{1}$.

It is worth noticing that $\Pi$ is necessarily a cycle of length $2$, \ie , $V_{t_{0}}\cdot l_1$ has  length $3$.
Note also that $\Pi$ consists of $2$ distinct points. Otherwise, $D$ is tangent to $l_1$ at $\Pi$,
so  this cycle $\Pi$ has to be supported at the  point of intersection
of $l_1$ with $l_2$ or $l_3$. But then, $Q$ does not lie on a unique cubic curve.

It follows from proposition \ref{prop:duality} that  $|{\I}_{P+q+\Pi }(5)|$ is $1$-irregular, provided it is complete.
This system clearly cannot contain any base curve, since $D$ is irreducible.
The potential extra base points are in the cycle $R=Q-\Pi$. If $R$ is composed of $7$
distinct points, it is clear that the irregularity of the system can be
at most $2$. Since the actual dimension of the system is at least $2$,
there can be at most $1$ extra base point, $r$.
Thus, the length $6$ cycle $R\setminus r$ is supported on one of the lines $l_2$ or $l_3$.
Since $D$ is irreducible, we get a contradiction.
If $R$ is not reduced, it can at most have one non reduced point at $l_2\cap l_3$.
Since the curve $D$ is irreducible, the linear system  $|{\I}_{P+q+\Pi }(5)|$ imposes
also a tangency condition at $r$ along the line $l_2$ (or $l_3$).
\qed

\end{pf}

\begin{rem} \label{rem:proj}It follows from the previous proof, that $S_5$
has no quadrisecant lines, if it is smooth.
For the same reason, if $S_5$ has excess of trisecant lines,
  a general trisecant line passing by a generic point
is not a quadrisecant line.
\end{rem}

\section{Applications of the existence theorem}
The existence theorem has an amusing consequence for the principal component
of the subvariety $W_{18}(5)$ of the Hilbert scheme $Hilb_{18}({\P}^2)$.

\subsection{The geometry of the generic point of  $W_{18}[5]$}

The space of length $18$ groups of points special in degree $5$, $W_{18} [5]$, 
is known to be irreducible of expected dimension $32$
\cite[theorem 3.2.1]{Co}, so it coincides with its principal component.
We shall show that its generic point has the same numerical
character as the projection by a generic trisecant line of the  generic White surface in ${\P}^5$.

\begin{cor} The generic point of $W_{18}[5]$ corresponds to a smooth uniform $18$-tuple of points, of numerical character
\[
\chi =(7,6,5,5,5);\]

therefore, it doesn't lie on any curve of degree less or equal to $4$.
\end{cor}

\begin{pf}  Let $\chi^{W} =(5+\epsilon_0 ,\cdots, 5+\epsilon_{4})$ be  the character of the projection of a generic White surface $S_5$
by a generic $3$-secant line. We have  $\sum \epsilon_{i} =3$. Since it corresponds to a $18$-tuple of points special in degree $5$,
$n_{0} \geq 7$.
If $n_0> 7$ or $n_i>6$ for $i>0$, we find $h^{1}(\I_{18}(5))\geq 2$. Thus, we have $\chi^{W} = (7,6,5,5,5)$.
Since $\chi^{W}$ is a uniform character, we find $dim( \chi )=32=dim(W_{18} [5])$. Furthermore, we have $H_{\chi^{W}}\subset W_{18}[5]$.
\qed
\end{pf}

Coppo's bound \cite{Co} shows, in this particular case, that the generic point of $W_{18}[5]$ corresponds to a smooth uniform $18$-tuple of points
not lying on any conic.

\subsection{ The number of $3$-secant lines passing through a generic point of $S$}

In 1882, H.~Krey \cite[p.505 , for $n=5$]{K} showed, using techniques  combining excess intersection theory and correspondence methods,
that, in degree $5$, there are $6$ associated pairs to $Z$, for $Z$ a generic $16$-tuple of points in the plane.
In this section only, a trisecant line to $S_5$ through a generic point $p=\Phi (q)$ is a line $|{\I}_{P+q+\Pi}(5)|^{\vee}$, where $\Pi$ is an associated
pair to ${\I}_{P}(5)$ at $q$. 
That is to say, we allow improper trisecant lines (\ie $\Phi (q)+\Phi (\pi_1)+\Phi (\pi_2)$ does not consist of distinct points on $S_5$). 
We prove that Krey's result still holds, if $Z$  corresponds  to the generic projection of any White surface $S_5$ from a generic point on it,
assuming that $S_5$ has a finite  number of trisecant lines passing by that point.

\begin{rem} In the case of a general polygonal surface, this  seems to be in apparent contradiction
with T.Dobler's result \cite[ proposition 3.17]{D}, which  shows that there are no trisecant lines at all.
But, the two  results agree. Since a polygonal surface in ${\P}^{5}$ has $6$ singular $4$-fold points,
the $6$  trisecant lines through a generic point $q$ that we construct in the next theorem are simply
the $6$ bisecant lines joining $q$ to a singular point of the surface. They are not counted by Dobler, for they are improper trisecant lines.
\end{rem}

Notice that  all trisecant lines to $S_5$ passing through a generic  point are obtained by the construction of theorem \ref{thm:1}.

\begin{thm} \label{thm:2}
Assume that, through a generic point $p:=\Phi(q)$ of a White surface $S_5$ in ${\P}^5$, there passes only a finite number of trisecant lines.
Then, through  $p$, there pass exactly $6$ trisecant lines counted with multiplicity.
\end{thm}

\begin{pf}
Let $D$ and $D^{\prime}$ be two generic quintic curves in $|{\I}_{P+q}(5)|$. Denote by $(D_\lambda)$ the pencil $<D,D^{\prime}>$.
From the construction of lemma \ref{lem:residu}, the cycle  that $D_\lambda$ induces on $D$, $R:=D\cdot D_{\lambda}$, lies on a unique cubic curve $E$.
If $D^{\prime}$ varies in $|{\I}_{P+q}(5)|$, the cubic curves $E$ vary in a linear system of
dimension at least $3$, with no multiple base  points. By Bertini's theorem,  for a generic choice of $D^{\prime}$, $E$ is smooth.

We show that no points of $P+q$ lie  on $E$, \ie , the pencil $D_{\lambda}$ induces a base point free $g_{6}^{1}$  on $E$.
This follows from the fact that $P+q$ imposes independent conditions on quartics.

\begin{lem}
The pencil $(D_{\lambda})$ cuts out on $E$ a base point free $g_{6}^{1}$.
\end{lem}
\begin{pf}
Indeed, suppose to the contrary, that the series $V_{\lambda}=(D_{\lambda}\cdot E)\setminus R$ has a base of length $l$.
First notice that $l<6$. Otherwise, applying the classical residuation theorem to $D\cdot E$ and $D\cdot D_{\lambda}$,
we deduce that $10$ of the points of $P+q$ lie on the same conic.
By genericity assumption, one can assume that 10 points of $P$ lie on the same conic;
so, $P$ lies on a quartic, giving a contradiction.

Let  $V_b$ denote the base of the series $V_{\lambda}$. Note that $V\setminus V_b$ lies on some conic $C_2^{\lambda}$.
Let us apply the residuation theorem with $m=5$, $n=3$, $n_1=5$ and $n_2=2$; we get

\[
D\cdot E=R+V_{b}+V^{\prime}\,\,\text{and}\,\, D\cdot D_{\lambda}=P+q+R=V_{b}+R+W
\]
\[
D\cdot C_2=V^{\prime}+T.
\]

We find $W+T=D\cdot C_{4}$, for some quartic $C_4$.
Suppose that $l=1$; then, the system of quartics passing through $W$ is empty, since $P+q$ imposes independent conditions.
Thus, $l\geq 2$. Notice that there is a pencil of conics $C_{2}^{\zeta}$ passing through $V^{\prime}$.
This pencil induces a linear series $T_{\zeta}$ on $D$ and a family of quartic curves $C_4^{\zeta}$, resolving the residuation theorem.

Suppose now that $l=2$; there is a unique quartic passing through $W$, so $T_{\zeta}$ is fixed  on $D$.
Thus, the pencil of conics has $10$ fixed points containing $9$ aligned points.
This line is then a fixed component of $D$, so $l\geq 3$.

Suppose that $l\geq 3$; we have, at least, a $(l-1)$-dimensional system of conics,
passing through $V^{\prime}$ and cutting out the series $T_{\zeta}$.
Since $P+q$ imposes independent conditions on quartics, this series is cut out by the $(l-2)$-dimensional system of quartic curves
passing through $W$. This leads to a contradiction, since $D$ cannot be a component of any of those systems.

Thus, the pencil $(D_{\lambda})$ induces a base point  free $g_{6}^{1}$ on $E$.\qed
\end{pf}

Let $L_{0}$ be the hyperplane section divisor on $E$ and  $L$  the divisor associated to this $g_{6}^{1}$.
We have an obvious morphism
\[
H^{0}(E,{\O}_{E}(L_{0}))\times H^{0}(E,{\O}_{E}(L_{1})) \xrightarrow{} H^{0}(E,{O}_{E}(L)),
\]
where $L_{1}=L-L_{0}$.
Its projectivization, ${\P}^{2}\times {\P}^{2}\xrightarrow {\phi} {\P}^{5}$, factors through the Segre embedding
of ${\P}^{2}\times {\P}^ {2}$ in ${\P}^{8}$ and a regular projection.
Its image is an hypersurface $\mathcal{A}$  of degree $6$ in ${\P}^{5}$,
  which represents the divisors of $|L|$ containing an aligned subscheme of length $3$.
Thus, the intersection of $\mathcal{A}$ with the line $D_{\lambda}$ parameterizes trisecant lines  to $S_5$
passing through $q$.
\qed
\end{pf}

\section{The  projection of $S_{5}$ by a generic trisecant}
In this section, we present Gambier's argument to show the finiteness of trisecant lines through a generic point of a White surface
\cite[p. 253-256]{G}. Gambier assumes implicitly that the triple curve $\gamma$ is irreducible; we show how to fill this gap.

\begin{thm}[Gambier] \label{thm:3}
The rational surface $\Sigma $, projection of $S_{5}$ from a  generic trisecant line, has degree $7$ and sectional genus $6$.
  On $\Sigma$, the double locus of the projection is in fact a triple  locus. Moreover, the curve $\gamma$ is a twisted cubic.
\end{thm}

\begin{pf}
Let $q$ be a generic point of ${\P}^2$ and $\Pi$ an associated pair to $|{\I}_{P}(5)|$ at $q$.
Denote by $Z$ the cycle $P+q+\Pi$ and by $\tilde{\P}\xrightarrow{\pi} {\P}^2$  the blowing up of ${\P}^2$ at $P+q+\Pi$.
We denote by $E$ the exceptional divisor of $\pi$.

Using Riemann-Roch's theorem, we find $\deg (\Sigma )=25-18=7$ and $\pi (\Sigma )=6$.
Only remains to prove that a generic hyperplane section of  $\Sigma $ contains $3$ triple points of the projection to ${\P}^{3}$.

A hyperplane section of $\Sigma $ has degree $7$, arithmetical genus $15$ and geometrical genus $6$, 
so it has either $9$ double points or $3$ triple
  points.
We shall construct these triple  points.
First, to construct a generic hyperplane section $h$ of $\Sigma$,  pick a generic point $s$ in ${\P}^{2}$.
The hyperplane section $h$  then corresponds  to a curve $D\in \left|{\I}_{Z+s}(5)\right|$.
According to Bertini's theorem, we may choose $D$ smooth, hence irreducible.
We write $\left|{\I}_{Z}(5)\right|=<D,D_1,\cdots D_3>$. By assumption, any $Y\in \chi (D,|{\I}_{Z}(5)|)$ lies on a unique conic curve.
An argument similar to the proof of proposition \ref{prop:19} shows that the strict transforms by $\pi$ of these conics 
vary in a linear system $\mathcal N$ on $\tilde{\P}$, as  $Y$ varies in the linear series cut out by $D^{\prime}\in <D_1,\cdots D_3>$. 
We find, as before, $dim(\mathcal{N})=dim(<D_1,\cdots ,D_3>)=2$.

Denote by $\tilde{D}$ the strict transform of $D$ by $\pi$ and consider the series $\mathcal{N}|_{\tilde D}$ cut on $\tilde D$ by $\mathcal{N}$.
Let  $\mathcal{L}$ denote the invertible sheaf of ${\O}_{\tilde{\P}}$-modules $\pi^{*}({\O}_{{\P}^2}(5))-E$.
By construction, we have $dim H^{0}(\mathcal{L}|_{\tilde D})=dim(\mathcal{N}|_{\tilde D})+1$.
Tensoring by $\mathcal{L}$ the exact sequence
\[
0\xrightarrow{} \mathcal{L}^{-1}\xrightarrow{} {\O}_{\tilde{\P}}\xrightarrow{} {\O}_{\tilde D}\xrightarrow{} 0 ,
\]
  we find
\[
\begin{aligned}
dim H^{0}(\mathcal{L}|_D)&=dim H^0(\mathcal{L})-1&=3\\
dim H^{1}(\mathcal{N}|_D)&=dim H^1(\mathcal{L})&=1.
\end{aligned}
\]

Thus, $dim H^0 (\mathcal{N}|_{\tilde{D}})=dim H^0(\mathcal{L}|_{\tilde D})=3$ and, by Riemann-Roch theorem, $\mathcal{L}|_{\tilde D}$ is a series of degree $7$.
The free part of $\mathcal{L}|_{\tilde D}$ is then a subseries of $\mathcal{N}|_{\tilde D}$, of the same dimension, $2$.
The two series are therefore equal. For any conic $C_2$, we have $D\cdot C_2=10$. So, $\mathcal{N}|_{\tilde D}$ has exactly $3$ base points
$A,B$ and $C$, distinct or infinitely near.
Since $\mathcal{N}$ is $2$-dimensional, the conics of this system cannot have a fixed line component,
so the support of the points $A,B,C$ consist of at least $2$ distinct points, say $A$ and $B$.
For the same reason, the points $A,B,C$, if distinct, are not aligned.

We write $D\cdot L= A+B+\lambda +\lambda^{\prime}+\lambda^{\prime\prime}$,
where $\lambda ,\lambda^{\prime}$ and $\lambda^{\prime\prime}$ are possibly
infinitely near.

\begin{lem} Let $L$ be the line joining $A$ and $B$;  the three points $\lambda ,
\lambda^{\prime}$ and $\lambda^{\prime\prime}$ distinct or infinitely near are mapped by  $\mid {\I}_{Z}(d)\mid $ to a triple point  of  $\Sigma$.

\end{lem}

\begin{pf} We only have to show that $\mid {\I}_{Z+\lambda}(5)\mid \subseteq \mid {\I}_{Z+\lambda +\lambda^{\prime}+\lambda^{\prime\prime}}(5)\mid$.
Let $r$ be any point of the plane. If  $\mid {\I}_{Z+r}\mid $ contains a pencil,  it is  complete, since a linear space of conics with fixed points, 
is at most $1$-dimensional.

Suppose that $\lambda$ is an ordinary point of the plane. Let $D^{\prime}$
be a generic curve of $|{\I}_{Z+\lambda}(5)|$. By the previous remark,  $D^{\prime}$  does not contain any of the points $A$, $B$ or $C$.
The cycle $Y:=(D^{\prime}\cdot D)\setminus Z$ lies on a unique conic $C_2$.
Since $C$ belongs to ${\mathcal N}$, the strict transform of this conic $C_2$ passes through $A,B$ and $C$.
Thus, the conic curve $C_2$ is the union of two lines $L$ and $l$, whose strict transforms pass through $C$.
By unicity of the conic $C_2$, two points of $Y$ must lie on $L$.
Thus, $\mid {\I}_{Z+\lambda}(5)\mid = \mid {\I}_{Z+\lambda +\lambda^{\prime}+\lambda^{\prime\prime}}(5)\mid$.
\qed
\end{pf}

By taking for $L$ the associated line and by working in the blow-up of ${\P}^2$ at $B$, the same proof works, if $C$ is infinitely near to $B$.
\qed
\end{pf}

\begin{rem} \hfil
\begin{enumerate}
\item Using Riemann-Roch's theorem on ${\mathcal N}\mid_{\tilde D^{\prime}}$, it is not hard to see that $A,B$ and $C$
are in fact the base points of $\mathcal N$.

\item The surfaces of ${\P}^{3}$  of degree $2n+3$,  passing $n$ times through a given twisted cubic,
have been studied by G\' erard\cite{Ge}.
For  $n=2$, we obtain another geometric construction of White surfaces, or more precisely, of their projections
to ${\P}^{3}$ from a trisecant line.

\end{enumerate}
\end{rem}

We have the following improvement of Dobler's result:

\begin{thm} The only White surface with a $4$-dimensional trisecant line locus is a Segre polygonal White surface.

\end{thm}

\begin{pf}  Let $\tilde S \xrightarrow{\pi} {\P}^{2}$ denote  the blow-up  of ${\P}^2$ at $Z$ for  $Z=P+q+\Pi$.
By remark \ref{rem:proj}, we  may assume that a generic trisecant line is not a quadrisecant line.
The only hypothesis of proposition \ref{prop:proj} that we still need to check, is that no plane curves are contracted by the projection from a
trisecant line passing through a generic point of $S_5$.

\begin{lem}If the White surface $S_5$ has an excess of trisecant lines, its projection to ${\P}^3$ from a generic trisecant line
does not contract any plane curve.
\end{lem}

\begin{pf}
Let $\phi_{S_5}$ (resp. $\phi_{\Sigma}$) denote the rational map associated to the linear system $|{\I}_{P}(5)|$ (resp. $|{\I}_{Z}(5)|$).
Let  $l:=<\phi_{S_5}(q),\phi_{S_5}(\pi_1), \phi_{S_5}(\pi_2)>$. Let $\pi_l$ denote the projection map from the trisecant line $l$,
and $\Sigma$ the image of $S_5$ by $\pi_{l}$.  The surface $\Sigma$ is the image of ${\P}^2$ by the rational map $\phi_{\Sigma}$.
 
 Assume, to the contrary, that there exists $\overline{C}$, a curve on $S_5$ which is contracted to a point $x$ of $\Sigma$
by the projection $\pi_{l}$. The curve $\overline{C}$ lies in the $2$-plane $<l,x>$.

\textbf{(A)} We may therefore assume that $\overline{C}$ is the $1$-dimensional part of the intersection $S_5\cap <l,x>$. 

We have $deg(\overline{C})\leq length(l\cap S_5)=3$.
The rational map $\phi_{S_5}$ is birational; let $C$ denote the plane curve $\phi_{\Sigma}^{-1}(\overline{C})$.
Then, the class of its strict transform $\tilde C$  by $\pi$, is given by 

\begin{equation}
[\tilde C] \cdot (5\pi^{*} L-\sum_{z\in Z} E_z)=0, \,\, \label{Eq:e}
\end{equation}
in $Pic(\tilde{S})$.
Furthermore, if $r_1$ and $r_2$ are two points of the curve $C$,
we have $|{\I}_{Z\cup\{r_1\}}(5)|=|{\I}_{Z\cup \{ r_2\}}(5)|$, since the points $r_1$ and $r_2$ are both mapped to
\[
x:=\phi_{\Sigma}(C)=\pi_{l}(\overline{C})\in\Sigma. 
\]
So, the curve $C$ must be a fixed component of the linear system of quintic curves passing through $r_1$ and $Z$. 
From this follows that the degree of $C$ must be strictly less than $5$.

We find  $[\tilde C]=d\pi^{*}L -\sum_{z\in Z^{\prime}} E_z$, where $d<5$ and $Z^{\prime}$ is a subset of $Z$.

From equation (\ref{Eq:e}) we get that the cardinality of $Z^{\prime}$ equals $5deg(C)$.
We denote by ${P}^{\prime}$ the subset of $P$ contained in $Z^{\prime}$. 
We get $deg(\overline{C})=5d-length(P^{\prime})$.
We conclude with a case by case study.

The case $d=4$ cannot happen since $Z$ has only $18$ points.
If $d=3$, then $C$ passes through at least $12$ points of $P$. Any $11$ points in ${P}^{\prime}$ do not lie on any conic,
for then $P$ would lie on a quartic curve. Hence, there is a unique cubic through those $11$ points.
We can assume that $q$ doesn't belong to that cubic, since $q$ is a generic point of the plane.
Therefore, $Z^{\prime}$ contains at least $13$ points of $P$, so that $P$ lies on a quartic curve.

Therefore we have $d\leq 2$.

If $d=1$, then the line $C$ passes through $5$ of the base points $Z$; so ${P}^{\prime}$ consists of at least $2$ points.
By genericity assumption, we can assume that $q$ doesn't lie on this line through ${P}^{\prime}$.
Thus, either $\overline{C}$ is a conic and $C$ is a line through $3$ points of $P$ and the associated pair $\{\pi_1, \pi_{2}\}$ to $q$
or $\overline{C}$ is a line and $C$ is a line through $4$ points of $P$ and only one point of the associated pair to $q$.

Suppose that $d=2$; from  equation (\ref{Eq:e}), we deduce that $C$ is a conic passing through $10$ points of $Z$,
among which $7$ at least belong to $P$. Since  no $6$ points of $P$ can be aligned, unless $P$ lies on a quartic curve,
we deduce that this conic curve $C$ through ${P}^{\prime}$ is unique. By genericity assumption, we may again assume that $q$ does not lie on $C$.
Then, either $\overline{C}$ is a conic and $C$ is a conic through $8$ points of $P$ and the associated pair to $q$
or $\overline{C}$ is a line and $C$ is a conic through $9$ points of $P$ and a single point of the associated pair to $q$.

Therefore, $\overline{C}$ is either a conic or a line not passing through $\phi_{S_5}(q)$.
Suppose that for $q$ generic in the plane and a generic trisecant line to $S_5$ through $q$,
there is a conic $\overline{C}$ contracted by $l$ to a point on $\Sigma$.
Then, $C$ is either a line through $3$ points of $P$ or a conic through $8$ points of $P$. Since there is a finite number of such curves,
as $q$ varies on $S_5$ there is a finite number of conics $\overline{C}$. 
From assumption \textbf{(A)}, $S_5$ is the union of a finite number of $1$-dimensional schemes. This gives an obvious contradiction.

Suppose now that for $q$ generic in the plane, and $l$ a generic trisecant line to $S_5$, there is a $2$-plane containing l meeting
$S_5$ along a $1$-dimensional locus, $\overline{C}(q,l)$ and  a finite number of points, such that $\overline{C}(q,l)$ is a line (not passing through $q$).
The curve $\overline{C}(q,l)$, corresponds in the plane to a curve $C(q,l)$, which is either a line through $4$ of the base points $P$ or a conic through $9$
of the base points $P$. Since there is a finite number of such curves, there is a finite number of curves $\overline{C}(q,l)$ as $q$ and $l$ vary.
Let us fix $q$ generic on $S_5$, by assumption there are infinitely many trisecant lines to $S_5$ through $q$. Therefore, there exists a trisecant line $l_0$
to $S_5$ through $q$ such that $\overline{C}(q,l_{0})$  meets an infinite number of trisecant lines to $S_5$ through $q$.
A line joining a generic point $\phi_{S_5}(\pi_1)$ of $\overline{C}(q,l_0)$ to $q$ is thus a trisecant line to $S_5$, so it meets $S_5$ at another point.
This third point, $\phi_{S_5}(\pi_2)$, varies in a curve as  $\phi_{S_5}(\pi_1)$ varies on $\overline{C}(q,l_0)$. From assumption \textbf{(A)},
 we find $deg(\overline{C}(q,l_0 ))>1$. This gives a contradiction.
\qed
\end{pf}

The class of the triple locus of the projection in $Pic({\tilde S})$ is therefore  
\[
9H^{*}-2\sum_{i=1}^{18}E_{i},
\] 
where $H^*$ is
the pull back to $\tilde {\P}$ of the divisor of lines in ${\P}^{2}$ and $E_i$ is the  exceptional divisor over the point $p_i$.
The curve $\gamma$ is therefore the birational transform of  a plane curve of degree $9$ passing twice through $(P+q+\Pi)$, which we still denote by $\gamma$.

  \begin{lem} If $S_5$ has a $4$-dimensional trisecant locus,  the curve $\gamma$ contains a line joining the points of $\Pi$,
some associated pair to $P+q$ in degree $5$.
\end{lem}

\begin{pf}
Suppose that $S_5$ has an excess of trisecant lines. 
Consider two pairs, $\Pi$ and $\Pi^ {\prime}$, associated to $|{\I}_{P}(5)|$ at $q$.
The points of $\Pi^{\prime}$ are mapped to the same point of $\Sigma$, so  the linear system $|{\I}_{Z+{\Pi}^{\prime}}(5)|$  is not complete.
By construction, any extra base point,  $\omega$, of $|{\I}_{Z+{\Pi}^{\prime}}(5)|$  lies on the line $<{\Pi}^{\prime}>$.
Since we can exchange $\Pi$ and ${\Pi}^{\prime}$, an extra base point, $\omega$, must lie at the intersection of $<\Pi >$ and $<{\Pi}^ {\prime}>$.
Thus, $S_{5}$ has only a finite number of trisecant  lines passing through our generic point, unless $<\Pi >$ is contained in $\gamma$.
\end{pf}

Assume that by a generic point of $S_{5}$ there pass  infinitely many trisecant lines. Then this is true at any point of the surface.
Since we can exchange the roles of $q,{\pi}_{1}$ and ${\pi}_ {2}$, the curve $\gamma$ is  the union of $3$ lines with a sextic curve,
contracted to a single point on ${\Sigma}$. Thus, $P$, the double points of $\gamma$  belong to this sextic, which is thus the union of $6$ lines,
meeting two by two at $P$.
Therefore, the White surface $S_5$, we started with, was a polygonal surface.
According to Dobler's thesis \cite[proposition 3.17]{D}, it is of Segre type.
\qed
\end{pf}

\begin{ack} This work grew out of precious geometric discussions I had with Professor Henry Pinkham,
during my Ph.D. studies at Columbia University and further stays in New York; I thank him for his kind support, during these years.
\end{ack}

\end{document}